\documentclass[journal]{IEEEtran}
\usepackage{cite}
\usepackage[pdftex]{graphicx}
\usepackage{amsmath}
\usepackage{algorithmic}
\usepackage{array}
\usepackage[utf8]{inputenc}
\usepackage[english]{babel}
\usepackage{amssymb}
\usepackage{algorithm2e}
\usepackage{gensymb}
\usepackage{amssymb}
 
\newtheorem{remark}{Remark}

\begin{document}

\title{Resilient Unit Commitment for Day-ahead Market Considering Probabilistic Impacts of Hurricanes}

\author{Tianyang~Zhao,~\IEEEmembership{Member,~IEEE,}
        Huajun~Zhang,~\IEEEmembership{Student Member,~IEEE,}
        Xiaochuan~Liu,~\IEEEmembership{Student Member,~IEEE,}
        Shuhan~Yao,~\IEEEmembership{Student Member,~IEEE,}
        and~Peng~Wang,~\IEEEmembership{Fellow,~IEEE}
\thanks{T. Zhao was with the Energy Research Institute, Nanyang Technological University, Singapore, e-mail: zhaoty@ntu.edu.sg.}
\thanks{H. Zhang, X. Liu, S. Yao and P. Wang were with the school
	of Electrical and Electronic Engineering, Nanyang Technological University, Singapore, e-mail: epwang@ntu.edu.sg}
}

\markboth{IEEE TRANSACTIONS ON POWER SYSTEMS}%
{Shell \MakeLowercase{\textit{et al.}}: Bare Demo of IEEEtran.cls for IEEE Journals}

\maketitle

\begin{abstract}
In the face of extreme events, e.g., hurricanes, the transmission systems, especially the transmission lines, are affected across time and space. To mitigate these impacts on the day-ahead market from a probabilistic perspective, a resilient unit commitment (UC) problem is formulated as a two-stage robust optimization (RO) problem. In the first stage, the status, energy, and reserves of generators are pre-scheduled to minimize the operational cost, responding to the worst line failure scenario in the operating day. The failure, operation status, and repair of transmission lines are depicted by a novel robust uncertainty set with chance constraint considering the repair of failed lines. This chance constraint is reformulated to its deterministic equivalence. Using both load shedding and generation curtailment, a recourse problem is formulated in the second stage considering the time-varying transmission lines operation status. The formulated RO problem is solved using a column-and-constraint generation scheme. Simulations are conducted on IEEE-24 and two-area IEEE reliability test system-1996 under hurricanes and results verify the effectiveness of the proposed method on the conservation of uncertainty set, worst-case line failure scenario detection and repair preparedness.
\end{abstract}


\IEEEpeerreviewmaketitle

\section*{Nomenclature}
\subsection{Indexes and sets}
\begin{IEEEdescription}[\IEEEusemathlabelsep\IEEEsetlabelwidth{\((i,j)(i,j)(i,j)\)}]
\item [\(g \in \mathcal{G}\)] Generator set
\item [\(d \in \mathcal{D}\)] Demand set
\item [\(i,j \in \mathcal{N}\)] Bus set
\item [\(ij \in \mathcal{E}\)] Transmission line set
\item [\(l \in \beth_{ij}\)] Conductor set of line \(ij\)
\item [\(k \in \daleth\)] Tower set
\item [\(t \in \mathcal{T}\)] Scheduling time periods
\end{IEEEdescription}
\subsection{Constants}
\begin{IEEEdescription}[\IEEEusemathlabelsep\IEEEsetlabelwidth{\((i,j)(i,j)(i,j)\)}]
\item[${s}_{\daleth,k}^t$] Equivalent wind speed rate at tower $k$ [m/s]
\item[\(\mu_{\daleth,k},\sigma_{\daleth,k}\)] Design parameter for tower $k$
\item[\(\pi_{\daleth,k}^{t}\)] Failure probability of tower $k$
\item[${L}_{ij,l}$]  Line segment length [km]
\item[${s}_{\beth_{ij},l}^t$,${S}_{ij}$] Wind speed rate and design wind speed rate at the location of segment $l$ [m/s]
\item[${Rf}_{\beth_{ij},l}^t$,${RF}_{ij}$] Rainfall rate and design rainfall rate at the location of segment $l$ [mm/h]
\item[${a}_{ij,l}$,${b}_{ij,l}$,${c}_{ij,l}$] Segment parameters
\item[\(\lambda_{\beth_{ij},l}^{t}\)] Failure rate of segment $l$
\item[\(\pi_{\beth_{ij},l}^{t}\)] Failure probability of segment $l$
\item[\(\pi_{ij}^{t}\)] Failure probability of line $ij$
\item[\(c_{\text{start},g}\)] Start-up cost of generator \(g\) [\$]
\item[\(c_{\text{shut},g}\)] Shut-down cost of generator \(g\) [\$]
\item[\(a_{g},b_{g}\)] Fuel cost of generator \(g\) [\$/MWh]
\item[\(c_{\text{R},g},c_{\text{r},g}^{+},c_{\text{r},g}^{-}\)] Spinning, regulation up and down reserve cost of generator \(g\) [\$/MWh]
\item[\(VOLL\)] Value of load loss [\$/MWh]
\item[\(VOGC\)] Value of generation curtailment [\$/MWh]
\item[\(\text{UT}_{g},\text{DT}_{g}\)] Minimum up/down time of generator \(g\) [h]
\item[\(\text{UT}_{\text{r}},\text{DT}_{\text{r}}\)] Remaining up/down duration of generator \(g\) [h]
\item[\(\text{RT}_{ij}\)] Repair duration of line \(ij\) [h]
\item[\(\text{R}^{60+}_{g},\text{R}^{60-}_{g}\)] Ramp up/down rate of generator \(g\) [MW/h]
\item[\(\text{SU}_{g},\text{SD}_{g}\)] Start-up and shut-down ramp limitation of \(g\) [MW]
\item[\(\text{R}^{10+}_{g}\)] Maximum 10-min ramp up rate of generator \(g\) [MW/h]
\item[\(\text{R}^{5+}_{g},\text{R}^{5-}_{g}\)] Maximum 5-min ramp up/down rate of generator \(g\) [MW/h]
\item[\(P_{g}^{\text{min}},P_{g}^{\text{max}}\)] Minimum/Maximum output of generator \(g\) [MW]
\item[\(P_{d}^{t}\)] Power demand of load \(d\) [MW]
\item[\(B_{ij}\)] Susceptance of line \(ij\) [S]
\item[\(T\)] Scheduling periods [h]
\item[\(D\)] Operating day
\item[\(\Delta t\)] Time step [h]
\item[\(\delta_{\text{R}}\)] Spinning reserve requirement
\item[\(\delta_{\text{r}}^{+},\delta_{\text{r}}^{+}\)] Regulation up/down reserve requirement
\item[\(M\)] big-M
\item[\(K\)] Maximum number of line failures
\item[\(\Pi\)]  Threshold probability of line failures
\item[\(\textbf{c},\textbf{d}\)]  Coefficient vector for the first/second stage objective function
\item[\(\textbf{G},\textbf{E},\textbf{M}\)]  Linear matrix for the first-stage, second-stage and uncertainty in the second-stage optimization problem.
\end{IEEEdescription}
\subsection{Uncertain Variables}
\begin{IEEEdescription}[\IEEEusemathlabelsep\IEEEsetlabelwidth{\((i,j)(i,j)(i,j)\)}]
\item [\(I_{ij}^{t}\)] Binary variable, 1 if line \(ij\) is off-line, 0 otherwise
\item [\(\kappa_{ij}^{t}\)] Binary variable, 1 if line \(ij\) is repaired, 0 otherwise
\item [\(\chi_{ij}^{t}\)] Binary variable, 1 if line \(ij\) is failed, 0 otherwise
\item [\(\boldsymbol{\xi}\)] Uncertain variable vector
\end{IEEEdescription}
\subsection{First-stage Decision Variables}
\begin{IEEEdescription}[\IEEEusemathlabelsep\IEEEsetlabelwidth{\((i,j)(i,j)(i,j)\)}]
\item[\(\alpha_{g}^{t}\)] Start-up command of generator \(g\)
\item[\(\beta_{g}^{t}\)] Shut-down command of generator \(g\)
\item[\(u_{g}^{t}\)] Binary variable, 1 if generator \(g\) is on-line,0 otherwise
\item[\(P_{g}^{t}\)] Energy set-point of generator \(g\) [MW]
\item[\(R_{g}^{t}\)] Spinning reserve of generator \(g\) [MW]
\item[\(r_{g}^{+t},r_{g}^{-t}\)] Regulation up/down reserve of generator \(g\) [MW]
\item[\(\theta_{i}^{t}\)] Angle at bus of line \(i\) [$\degree$]
\item[\(P_{ij}^{t}\)] Power transfer on line \(ij\) [MW]
\item[\(Q_{\text{OR}}^{t}\)] Operation reserve [MW]
\item[\(\textbf{x}\)] First-stage decision variable vector
\end{IEEEdescription}
\subsection{Second-stage Decision Variables}
\begin{IEEEdescription}[\IEEEusemathlabelsep\IEEEsetlabelwidth{\((i,j)(i,j)(i,j)\)}]
\item[\(p_{g}^{t}\)] Active power output of generator \(g\) [MW]
\item[\(p_{\text{c},g}^{t}\)] Generation curtailment of generator \(g\) [MW]
\item[\(p_{ij}^{t}\)] Power transfer on line \(ij\) [MW]
\item[\(\gamma_{i}^{t}\)] Angle at bus of bus \(i\) [$\degree$]
\item[\(p_{d,j}^{t}\)] Load shedding of load \(d\) [MW]
\item[\(\textbf{y}\)] Second-stage decision variable vector
\end{IEEEdescription}

\section{Introduction}
\IEEEPARstart{H}{urricanes} can damage the transmission network assets, e.g., towers and conductors, affecting the reliability of power systems across time and space\cite{zhang2019spatial}. These impacts can be mitigated before, during and after the hurricanes\cite{wang2015research}, calling for the resilience management of power systems towards extreme events. Before the advent of one event, for the uncertain failures of components during this event, the schedulable resources should be well prepared to avoid both load shedding and generation curtailment\cite{zhao2017distributionally} within different operation processes, e.g., day-ahead market and real-time market. As the bridge between the available resources and real-time market, the day-ahead market should not only clear the market for the given energy and ancillary requirements but also provide guidance to market players in advance toward the uncertain contingencies\cite{ela2014evolution}.

Considering the mathematical characteristics of contingency uncertainties, the stochastic\cite{fernandez2016probabilistic,zheng2014stochastic,Dimitris2019Resilient,wang2016resilience,wang2018resilience}, robust\cite{wang2013two,guo2016contingency,hu2015robust,wang2016robust} and distributionally robust\cite{zhao2017distributionally} unit commitment (UC) problems are formulated to schedule the generators, transmission lines, and demand side resources to improve the reliability and operational efficiency. In the stochastic UC, the contingencies of generations and transmission lines have been modeled as probabilistic sets, e.g., possibility distribution functions (PDFs) with fixed parameters. These PDFs have formulated of UC problems under $N$-$k$ or $N$-1 security criteria regarding the loss of load
probability (LOLP) and expected energy not served
(EENS) in \cite{fernandez2016probabilistic,zheng2014stochastic}. Capturing the time-varying operation features, the time-varying PDFs have been introduced to improve the resilience of day-ahead operation in\cite{Dimitris2019Resilient}. Updating the probability under the real-time operation conditions, the real-time management of power systems, e.g., hourly ahead UC, has been formulated as a dynamic process in\cite{wang2016resilience,wang2018resilience}. 

Without using the probability explicitly, the robust UC problems are formulated with a polyhedral uncertainty set to depict the combination of contingencies. A two-stage robust UC is proposed to consider the uncertainty of $N-k$ contingencies of both generators and transmission lines in \cite{wang2013two}. To meet the $N$-1-1 contingency reliability criterion, a time dependent operation status uncertainty set for generators and transmission lines is proposed in \cite{guo2016contingency}. Aside from the uncertain generation contingencies, the interval uncertainty of loads is covered in \cite{hu2015robust}. The corrective actions have been integrated into the robust UC framework to enhance the resilience of energy systems in \cite{wang2016robust}.

Intersecting of probabilistic set and robust uncertainty set, the ambiguity set can improve the robustness of the probabilistic sets using distributionally robust optimization. A distributionally robust contingency constrained UC is proposed to manage the ambiguity of time-invariant failures probabilities in\cite{zhao2017distributionally}. A follow-up work is proposed to manage the distribution network reconfiguration towards random contingencies depicted by the ambiguity set in\cite{babaei2018distributionally}.

Under hurricanes, i.e., a time-varying track and intensity, the transmission networks, e.g., towers and conductors, might be destroyed at different time and location along with the evaluation of one hurricane\cite{zhang2019spatial,hurricane_report}. After failure, one line will stay off-line until the repair. The mode transition, i.e., on-line \(\xrightarrow{\text{Failure}}\) off-line \(\xrightarrow{\text{Repair}}\) on-line, asks for a time-varying uncertain model for lines. What is more, during the hurricane, the repair usually takes several hours or days\cite{subcommittee1979ieee}, indicating one line might remain off-line to the end of scheduling once failed at any time slot\cite{hurricane_report,van2015transmission}. The failed time slot may not always be the start of the scheduling periods\cite{hurricane_report}, and varies for different failed lines responding to the track of hurricanes. The state-of-art time-varying line failure models have considered the probability of failure\cite{Dimitris2019Resilient} and operation status\cite{Dimitris2019Resilient,wang2016resilience,wang2018resilience}, while the repair has not been considered.

With the deregulation of power systems, the day-ahead market is playing important roles in scheduling\cite{ela2014evolution,eu_market_report}, and works within extreme weather events\cite{noauthor_despite_nodate}. As a transmission network dependent\footnote{The change in the transmission lines might affect the power markets, e.g., line impedance, network typologies.} management strategy, the day-ahead market should be cleared to meet the energy and ancillary requirements while considering the uncertain contingencies. The cleared market results should provide appropriate notice and warning to market participants and other relevant entities\cite{hurricane_report}. Under some extreme line failures, the generators along the hurricane path might be forced to shut down\cite{sandi_report}. The day-ahead market should be extended to consider the uncertain line failures before the advent of hurricanes.

Oriented from the probabilistic impacts of hurricanes on transmission lines, a robust uncertainty set with chance constraints is proposed for the transmission lines with failure and repair. This uncertainty set is integrated into a two-stage robust optimization (RO) problem to minimize the total cost in the day-ahead market. In the first stage, the market is cleared under normal conditions, and the second stage reveals the impacts of the worst line failure scenario on load shedding and generation curtailment. The RO problem is solved using the column-and-constraint generation (C\&CG) scheme. The contribution of this paper can be summarized as follows: 1) a robust uncertainty set with chance constraints is proposed for the transmission line failures under hurricane, with failure and repair, 2) a two-stage robust UC is formulated for the day-ahead market using the uncertainty set to minimize the total cost under worst line failure scenario.

The rest of this paper is organized as follows: the resilience management for the day-ahead market and probabilistic line failures under hurricanes are introduced in Section II. The two-stage UC is formulated in Section III, with the proposed robust uncertainty set under chance constraints. The C\&CG scheme is presented in Section IV. Case studies are conducted in Section V, and conclusions are drawn on Section VI.

\section{Resilient Day-ahead Market Operation Considering Impacts of Hurricanes}
Considering the impacts of hurricanes, the system operator schedules the available resources in the day-ahead market to enhance the resilience of power systems. The impacts of hurricanes on transmission line failures are depicted as the PDFs in this section.
\subsection{Resilience Management of Day-ahead Market}
\begin{figure}[!h]
\includegraphics[width=1\linewidth]{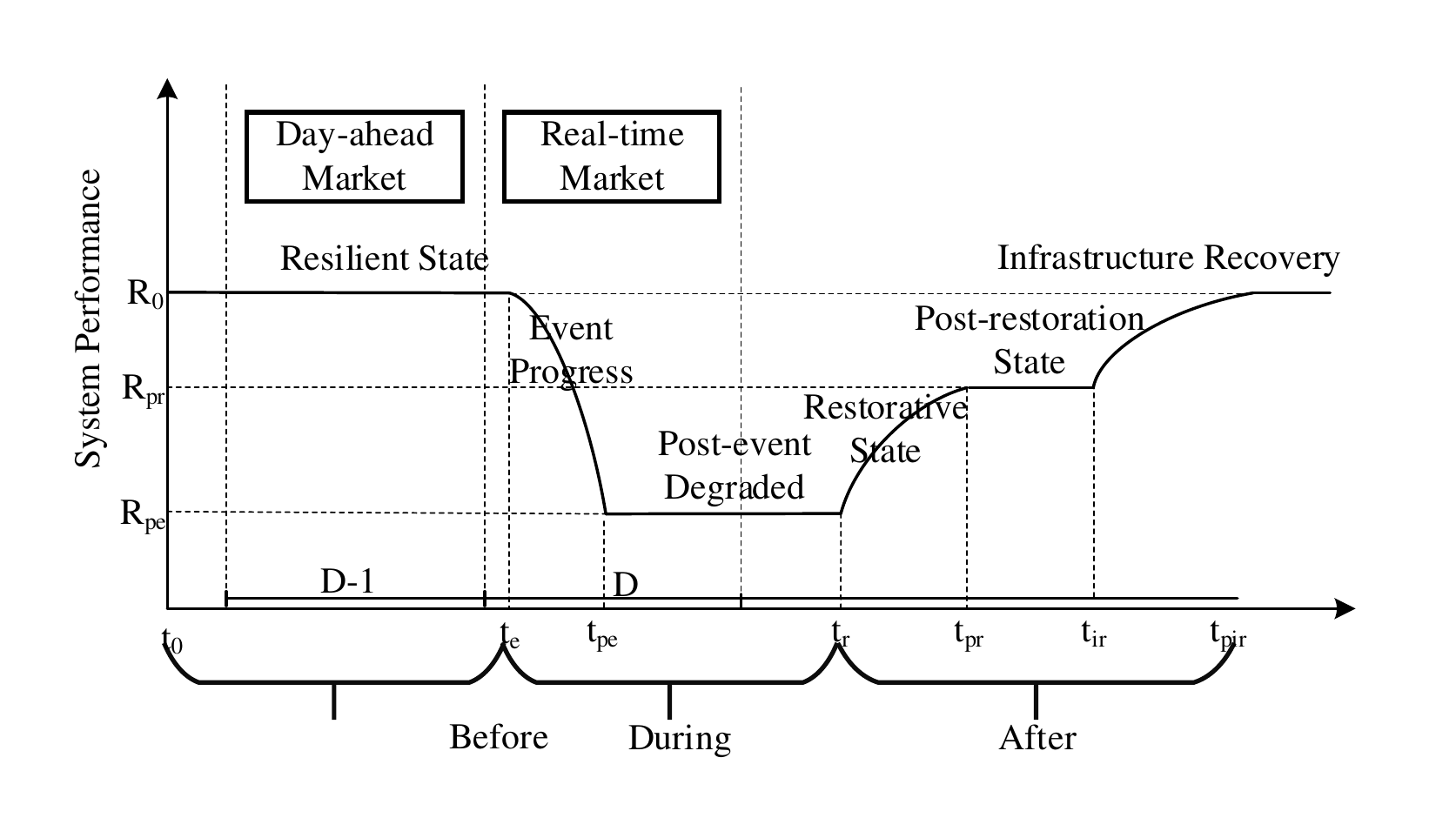}
\caption{Resilient management of power markets towards hurricanes.$[t_{\text{0}},t_{\text{e}}]$ is the pre-disturbance resilient state. $[t_{\text{e}},t_{\text{pe}}]$ is the event process. 
$[t_{\text{pe}},t_{\text{r}}]$ is the post-event degraded state.
$[t_{\text{r}},t_{\text{pr}}]$ is the restorative state.  $[t_{\text{pr}},t_{\text{ir}}]$ is the post-restoration state. $[t_{\text{ir}},t_{\text{pir}}]$ is the infrastructure recovery.}
\label{fig:resi_mag}
\end{figure}

Consider a transmission systems with a set of generators, loads, buses and transmission lines, denoted by \(\mathcal{G}\), \(\mathcal{D}\), \(\mathcal{N}\), and \(\mathcal{E}\), respectively. A hurricane is forecast to land across a given time horizon in the following day, i.e., \([t_{\text{e}},t_{\text{pe}}]\) in Fig.\ref{fig:resi_mag}. The operating day horizon is denoted by \(\mathcal{T}\), and discrete into equal time slot by time step \(\Delta t\). 

As shown in Fig.\ref{fig:resi_mag}, the resilience management of one power system is generally classified into three stages, i.e., before, during and after\cite{wang2015research,yao2018transportable}. Before the advent of this hurricane, i.e. $\text{resilient state} [t_{0},t_{\text{e}}]$\footnote{$t_{\text{e}}$ is assumed to belong to $\mathcal{T}$, while $t_{\text{pe}}$ is not necessarily within $\mathcal{T}$.}, the generators, i.e., \(g \in \mathcal{G}\), are scheduled to maximize the operational resilience or minimize the total operational cost. In the day-ahead market, after receiving the bids and offers from different market players, the power market is cleared responding to the load demand and ancillary requirements, considering the worst system performance, i.e., $\text{R}_{\text{pe}}$ in Fig.\ref{fig:resi_mag}, responding to probabilistic impacts of hurricanes on transmission lines within $\mathcal{T}$, as depicted in Section II.B. The market is cleared using the resilient UC model in Section III. After the clear of power market, the commitment, energy, reserves and worst scenario generation curtailment of generators are passed to the reliability assessment commitment\cite{zhang2019spatial,singhal_hedman_2019} or real-time market. In this paper, the day-ahead market is to be analyzed.
\subsection{Probabilistic Impacts of Hurricanes on Transmission Lines}
\begin{figure}[!h]
\centering
\includegraphics[width=0.9\linewidth]{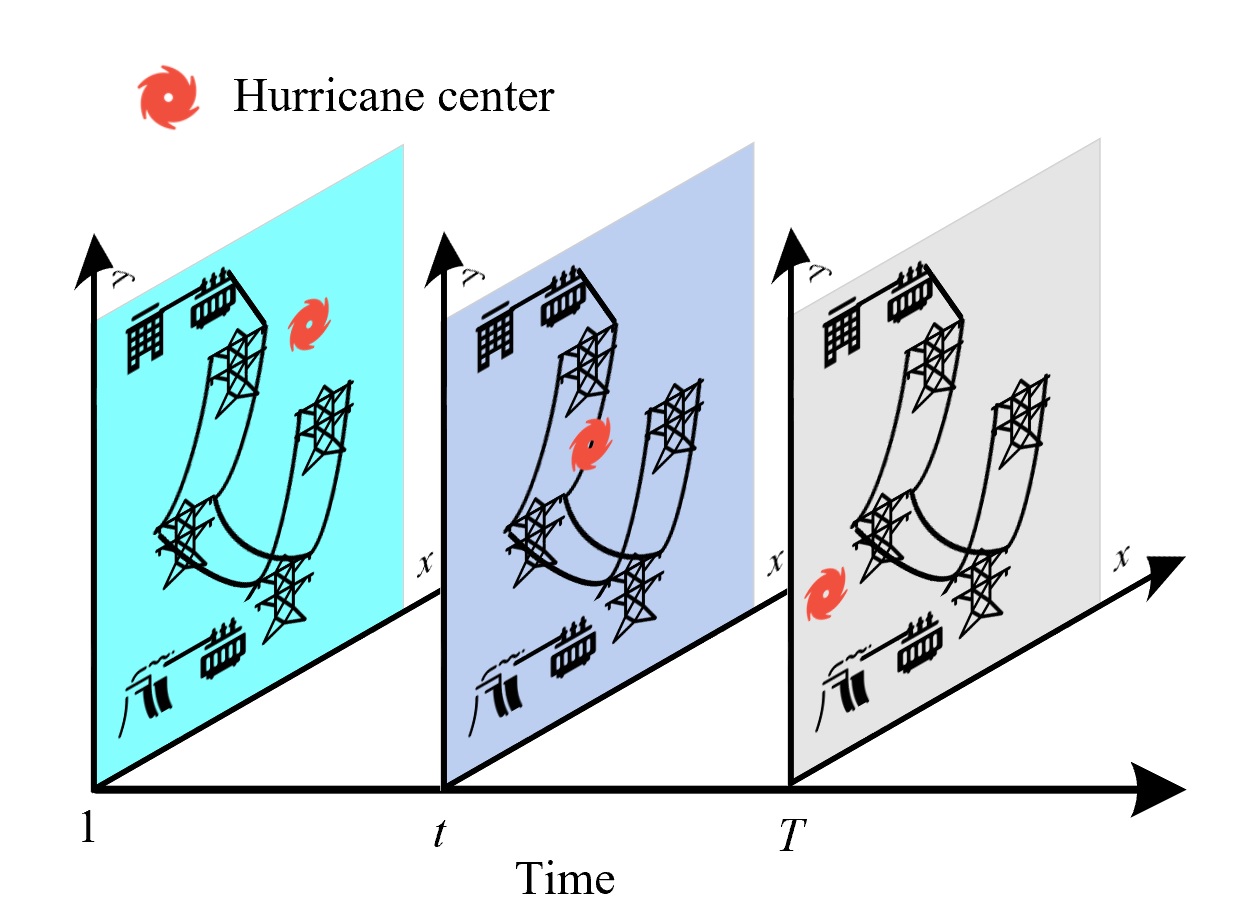}
\caption{Spatial and temporal impacts of hurricanes on transmission lines. The track of one hurricane is depicted by the movement of hurricane center across time and space,i.e.,$(x_{1},y_{1},1)\rightarrow (x_{t},y_{t},t) \rightarrow (x_{T},y_{T},T)$.}
\label{fig:hurricane}
\end{figure}
A probabilistic line failure model is formulated to quantify the impacts of hurricanes on transmission lines from both spatial and temporal perspectives. The hurricane is treated as a time-varying track and intensity across the given area, i.e., one power system in this paper, as shown in Fig.\ref{fig:hurricane}. Without considering the flood along with hurricane, the failure caused by the hurricanes is mainly transmission line related, e.g., the transmission towers and conductors\cite{hurricane_report}. 

In this sub-section, one analytical approach is adopted to quantify the impacts of hurricanes on transmission lines. For line $ij \in \mathcal{E}$, it consists of several towers and conductors, i.e., $\beth \cup \daleth$. The time-varying failure probability of tower $k$ is given as follows 
\begin{equation}
\begin{aligned}
\label{eq:failure_prob_tower}
\pi_{\daleth,k}^t = \int_{-\inf}^{{s}_{\daleth,k}^t}  \frac{1}{{\sigma_{\daleth,k} \sqrt {2\pi } }}e^{{{ - \left( {x - \mu_{\daleth,k} } \right)^2 } \mathord{\left/ {\vphantom {{ - \left( {x - \mu } \right)^2 } {2\sigma ^2 }}} \right. \kern-\nulldelimiterspace} {2\sigma_{\daleth,k} ^2 }}},\forall t,k
\end{aligned}
\end{equation}

For conductors, considering the line length is much larger than the dimension of the hurricane center, the conductor for line $ij$ is divided into an appropriate number of equal segments, denoted by set \(\daleth_{ij}\). The regression model\cite{Xiao2006c} is extended to take into account of rainfall impact in addition to strong wind. The time-varying failure rate of segment $l$ is expressed as follows
\begin{equation}
\begin{aligned}
\label{eq:failure_rate_cs}
{\lambda}_{\beth_{ij},l}^t={L}_{ij,l}\exp\left(\frac{{a}_{ij,l}{s}_{\beth_{ij},l}^t}{{S}_{ij,l}}+\frac{{b}_{ij,l}{Rf}_{\beth_{ij},l}^t}{{RF}_{ij}}+{c}_{ij,l}\right),\\
\ \forall t, ij, l
\end{aligned}
\end{equation}

According to the discrete time Markov process with constant failure rate at time slot $t$, the failure probability ${\pi}_{ij,k}^t$ of segment $k$ can be calculated as follows \cite{Liu2011f}
\begin{equation}
\begin{aligned}
\label{eq: failure_prob_cs}
\pi_{\beth_{ij},l}^t=(1-\pi_{\beth_{ij},l}^{t\text{\scriptsize{-}}1})(1-\exp(-\lambda_{\beth_{ij},l}^t\Delta t))+\pi_{\beth_{ij},l}^{t\text{\scriptsize{-}}1}, \forall t,ij,l
\end{aligned}
\end{equation}

When the transmission lines and towers fail independently and any failure of transmission segment or tower results in the failure of line $ij$, its failure probability is given as follows

\begin{equation}
\begin{aligned}
\label{eq:failure_prob_line}
\pi_{ij}^t= 1-\prod_{k \in ij}(1-\pi_{\daleth,k}^t)\prod_{l \in \beth_{ij}}(1-\pi_{\beth_{ij},l}^t), \forall t, ij
\end{aligned}
\end{equation}

In Eq.(\ref{eq:failure_prob_line}), the failure probability follows a Bernoulli distribution with parameter $\pi_{ij}^t$., i.e., $\chi_{ij}^{t} \sim \text{B}(1,\pi_{ij}^t)$.
\begin{remark}
In Eq.(\ref{eq:failure_prob_tower})-(\ref{eq:failure_rate_cs}), the equivalent wind speed rate and rainfall speed, i.e., ${s}_{\daleth_{ij},k}^t$, ${s}_{ij,k\beth_{ij},l}^t$ and ${Rf}_{\beth_{ij},l}^t$, depend on the distance between the segment $k,l$ and the hurricane center\cite{zhang2019spatial}. If the tower or segment is too far away from the hurricane center, the failure probabilities of tower and segments will not be affected by the hurricane.
\end{remark}

\section{Resilient Unit Commitment with Probabilistic Line Failures}
Considering the uncertainty of line failures, a two-stage RO problem is formulated in this section, where the market is cleared under normal conditions in the first stage and the worst system performance is assessed in the second stage.
\subsection{Objective function}
The objective function is to minimize the total operational cost, i.e., start-up, shut-down, fuel and reserve cost of generators, together with the worst load shedding and generation curtailment cost considering the operation status of transmission lines, as follows
\begin{equation}\label{eq:RO}
\begin{aligned}
\min_{\textbf{x} \in \textbf{X}} f(\textbf{x})+\max_{\boldsymbol{\xi} \in \mathcal{U}}[\mathcal{Q}(\textbf{x})]
\end{aligned}
\end{equation}
\begin{equation}\label{eq:obj}
\begin{aligned}
& f(\textbf{x}) = \textbf{c}^{\text{T}}\textbf{x} \\ 
& = \sum_{t \in \mathcal{T}}\sum_{g \in \mathcal{G}} \{ \underbrace{c_{\text{start},g}\alpha_{g}^{t} + c_{\text{shut},g}\beta_{g}^{t}}_{\text{start-up and shut-down cost}} \\ 
& + \underbrace{b_{g}u_{g}^{t} + [a_{g}P_{g}^{t}}_{\text{fuel cost}} + \underbrace{c_{\text{R},g}R_{g}^{t} + c_{\text{r},g}^{+}r_{g}^{+t} + c_{\text{r},g}^{-}r_{g}^{-t}}_{\text{reserve cost}}]\Delta t\}
\end{aligned}
\end{equation}
\begin{equation}\label{eq:second_stage}
\begin{aligned}
&\mathcal{Q}(\textbf{x}) = \min_{\textbf{y}\in \textbf{Y}} \textbf{d}^{\text{T}}\textbf{y} \\
& = VOLL\sum_{t \in \mathcal{T}}\sum_{d \in \mathcal{D}} p_{d}^{t}
+VOGC\sum_{t \in \mathcal{T}}\sum_{g \in \mathcal{G}} p_{\text{c},g}^{t} \Delta t \\
\end{aligned}
\end{equation}
where \(\textbf{x}:=\{\alpha_{g}^{t},\beta_{g}^{t},u_{g}^{t},P_{g}^{t},R_{g}^{t},r_{g}^{+t},r_{g}^{-t},P_{ij}^{t},\theta_{i}^{t},Q_{\text{OR}}^{t}\}\) is the first stage decision variable vector. \(\textbf{y}:=\{p_{g}^{t}, p_{d}^{t}, \gamma_{i}^{t}, p_{ij}^{t},p_{\text{c},g}^{t}\}\) is the second stage decision variable vector. \(\boldsymbol{\xi}:=\{I_{ij}^{t}, \chi_{ij}^{t}, \kappa_{ij}^{t}\}\) is the uncertain variable vector. \(\textbf{X}\) is the constraint set for the first-stage scheduling, depicted in Section III.B. \(\textbf{Y}:=\{\textbf{E}\textbf{y} \geq \textbf{h}-\textbf{G}\textbf{x}-\textbf{M}\boldsymbol{\xi}\}\) is the constraint set for real-time operation after the realization of transmission line operation status, depicted in Section III.C. The constraint set for uncertain variables, i.e., \(\mathcal{U}\) is illustrated in Section III.D.

\subsection{First-stage Constraint Set}
The first-stage constraint set, i.e., \(\textbf{X}\), is to meet the normal operation requirement, including the technical constraints of generators, transmission lines, power balance and reserve requirements are formulated as follows
\begin{equation}\label{eq:gen_status}
\alpha_{g}^{t}-\beta_{g}^{t} = u_{g}^{t}-u_{g}^{t-1}, \forall t, g
\end{equation}
\begin{equation}\label{eq:up_duration}
\sum_{q=t-\text{UT}_{g}+1}^{t} \alpha_{g}^{q} \leq u_{g}^{t}, \forall t \in \{\text{UT}_{g},...,T\}, g
\end{equation}
\begin{equation}\label{eq:down_duration}
\sum_{q=t-\text{DT}_{g}+1}^{t} \beta_{g}^{q} \leq 1-u_{g}^{t}, \forall t \in \{\text{DT}_{g},...,T\}, g 
\end{equation}
\begin{equation}\label{eq:initial_status}
u_{g}^{t} = u_{g}^{0}, \forall t \in \{\Delta t,...,\text{UT}_{\text{r}}+\text{DT}_{\text{r}}\},g
\end{equation}
\begin{equation}\label{eq:power_lower}
u_{g}^{t}P_{g}^{\text{min}} \leq P_{g}^{t} - r_{g}^{-t}  , \forall t,g 
\end{equation}
\begin{equation}\label{eq:power_upper}
P_{g}^{t} + R_{g}^{t} + r_{g}^{+t}\leq P_{g}^{\text{max}}u_{g}^{t}, \forall t,g
\end{equation}
\begin{equation}\label{eq:spinning}
0 \leq R_{g}^{t} \leq \text{R}_{g}^{10+}u_{g}^{t}, \forall t,g
\end{equation}
\begin{equation}\label{eq:regulation_up}
0 \leq r_{g}^{+t} \leq \text{R}_{g}^{5+}u_{g}^{t}, \forall t,g
\end{equation}
\begin{equation}\label{eq:regulation_down}
0 \leq r_{g}^{-t} \leq \text{R}_{g}^{5-}u_{g}^{t}, \forall t,g
\end{equation}
\begin{equation}\label{eq:ramp_up}
P_{g}^{t} - P_{g}^{t-1} \leq \text{R}_{g}^{60+}u_{g}^{t-1}\Delta t+\alpha_{g}^{t}\text{SU}_{g}, \forall t, g
\end{equation}
\begin{equation}\label{eq:ramp_down}
P_{g}^{t-1} - P_{g}^{t} \leq \text{R}_{g}^{60-}u_{g}^{t}\Delta t+\beta_{g}^{t}\text{SD}_{g}, \forall t, g
\end{equation}
\begin{equation}\label{eq:power_balance}
\begin{aligned}
\sum_{g \in \mathcal{G}_{j}} P_{g}^{t}+ \sum_{ij} P_{ij}^{t} - \sum_{ji}P_{ji}^{t} = \sum_{d \in \mathcal{D}_{j}} P_{d}^{t}, \forall t, j
\end{aligned}
\end{equation}
\begin{equation}\label{eq:branch}
\begin{aligned}
P_{ij}^{t}-B_{ij} \left(\theta_{i}^{t}-\theta_{j}^{t}\right)=0, \forall t, ij
\end{aligned}
\end{equation}
\begin{equation}\label{eq:power_flow}
P_{ij}^{\min} \leq P_{ij}^{t} \leq P_{ij}^{\max}, \forall t, ij
\end{equation}
\begin{equation}\label{eq:operation_reserve}
Q_{\text{OR}}^{t} \geq P_{g}^{t} + R_{g}^{t}, \forall t,g
\end{equation}
\begin{equation}\label{eq:spinning_reserve_requirement}
\sum_{g \in \mathcal{G}} R_{g}^{t} \geq \delta_{\text{R}} Q_{\text{OR}}^{t}, \forall t
\end{equation}
\begin{equation}\label{eq:re_up_requirement}
\sum_{g \in \mathcal{G}} r_{g}^{+t} \geq \delta_{\text{r}}^{+}\sum_{d \in \mathcal{D}} P_{d,j}^{t}, \forall t
\end{equation}
\begin{equation}\label{eq:re_down_requirement}
\sum_{g \in \mathcal{G}} r_{g}^{-t} \geq \delta_{\text{r}}^{-}\sum_{d \in \mathcal{D}} P_{d,j}^{t}, \forall t
\end{equation}

Eq.(\ref{eq:gen_status}) is the operation status transition of generators. Eq.(\ref{eq:up_duration})-(\ref{eq:initial_status}) are the minimal up/down time duration, and initial status constraint of generators\cite{morales2012tight}. Eq.(\ref{eq:power_lower})-(\ref{eq:regulation_down}) are the power capacity limitation with reserves\cite{morales2012tight,li2015flexible}. Eq.(\ref{eq:ramp_up})-(\ref{eq:ramp_down}) are the ramp up and down limitation \cite{li2015flexible}. The power balance of each bus is given in Eq.(\ref{eq:power_balance}). The power transmitted on each line is given in Eq.(\ref{eq:branch}), and is limited by Eq.(\ref{eq:power_flow}). Considering the uncertainty of load forecasting and possible failure of any generator, the spinning and regulation reserve requirements are given in Eq.(\ref{eq:operation_reserve})-(\ref{eq:re_down_requirement}).

\subsection{Second-stage Constraint Set}
The second stage constraint set, i.e., \(\textbf{Y}\), includes the real-time scheduling of generators, power flow, load shedding, generation curtailment after the realization of transmission line operation status, as follows 
\begin{equation}\label{eq:gen_sec}
\begin{aligned}
P_{g}^{t} - R_{g}^{t} \leq p_{g}^{t} \leq P_{g}^{t} + R_{g}^{t}, \forall t, g
\end{aligned}
\end{equation}
\begin{equation}\label{eq:ramp_sec}
\begin{aligned}
-\text{R}_{g}^{60-}\Delta t \leq p_{g}^{t}-p_{g}^{t-1} \leq \text{R}_{g}^{60+}\Delta t, \forall t, g
\end{aligned}
\end{equation}
\begin{equation}\label{eq:branch_sec}
\begin{aligned}
\left(I_{ij}^{t}-1\right)M \leq p_{ij}^{t}-B_{ij} \left(\gamma_{i}^{t}-\gamma_{j}^{t}\right) \leq \left(1-I_{ij}^{t}\right)M,\forall t, ij
\end{aligned}
\end{equation}
\begin{equation}\label{eq:power_flow_sec}
\begin{aligned}
-I_{ij}^{t} P_{ij}^{\max} \leq p_{ij}^{t} \leq I_{ij}^{t} P_{\text{ij}, k}^{\max}, \forall t, ij
\end{aligned}
\end{equation}
\begin{equation} \label{eq:load_shedding}
P_{d}^{t} \geq p_{d}^{t} \geq 0, \forall t, d
\end{equation}
\begin{equation} \label{eq:generation_curtailment}
P_{g}^{t} \geq p_{\text{c},g}^{t} \geq 0, \forall t, g
\end{equation}
\begin{equation}\label{eq:power_balance_sec}
\begin{aligned}
\sum_{g \in \mathcal{G}_{j}} (p_{g}^{t}-p_{\text{c},g}^{t})+\sum_{ij} p_{ij}^{t}-\sum_{ji} p_{ji}^{t}= \sum_{d \in \mathcal{D}_{j}} (P_{d}^{t}-p_{d}^{t}),  
\forall t, j 
\end{aligned}
\end{equation}

The real-time generator output constraints are given in Eq.(\ref{eq:gen_sec})-(\ref{eq:ramp_sec}). Considering the operation status of transmission lines, the power transmitted on each line is limited by Eq.(\ref{eq:branch_sec})-(\ref{eq:power_flow_sec}). For specific bus, the amount of load shedding is given in Eq.(\ref{eq:load_shedding}). The generation curtailment is limited by constraint (\ref{eq:generation_curtailment}). The real-time power balancing at each bus is depicted by Eq.(\ref{eq:power_balance_sec}). 

\subsection{Uncertainty Set}
\begin{figure}[!h]
\centering
\includegraphics[width=0.7\linewidth]{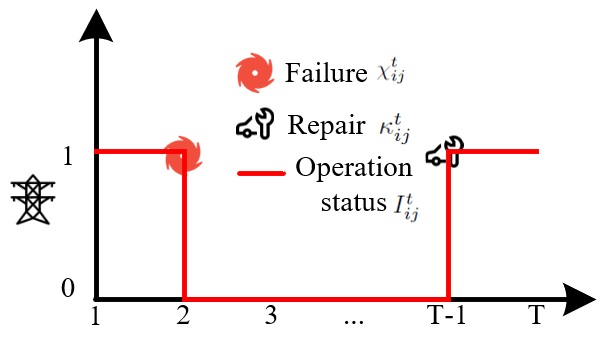}
\caption{A conceptual description on failure, operation and repair of one transmission line.}
\label{fig:failure}
\end{figure}

As shown in Section II.B and Section III.C, the failure $\chi_{ij}^t$ and operation status of line $I_{ij}^{t}, ij \in \mathcal{E}$ are given, respectively. Clearly, for any transmission line, its operation status depends on both failure and repair, as shown in Fig.\ref{fig:failure}. This relationship can be depicted by the following equation
\begin{equation}\label{eq:D1}
\begin{aligned}
\kappa_{ij}^{t}-\chi_{ij}^{t} = I_{ij}^{t} - I_{ij}^{t-1}, \forall t, ij
\end{aligned}
\end{equation}

During the scheduling periods, \(\mathcal{T}\), any transmission line $ij$ is assumed to be failed at most once when the repair is not considered, as follows 
\begin{equation}\label{eq:D2}
\begin{aligned}
\sum_{t \in \mathcal{T}} \chi_{ij}^{t} \leq 1, \forall ij
\end{aligned}
\end{equation}
\begin{equation}\label{eq:D21}
\begin{aligned}
\kappa_{ij}^{t}=0, \forall t, ij
\end{aligned}
\end{equation}

Considering the impacts of hurricanes as the line failure probability in Section II.B, the following constraint is proposed to detect which line is affected by the hurricane during which time slots via introducing a probabilistic threshold \(\Pi\)
\begin{equation}\label{eq:D3}
\begin{aligned}
\chi_{ij}^{t}(\pi_{ij}^{t}-\Pi) \geq 0,\forall t, ij
\end{aligned}
\end{equation}

In preparation of one hurricane, at most $K$ transmission lines are considered to be off-line within one time slot, as follows
\begin{equation}\label{eq:D4}
\begin{aligned}
\sum_{ij\in \mathcal{E}} I_{ij}^{t} \geq |\mathcal{E}|-K, \forall t
\end{aligned}
\end{equation}

In summary, \(\mathcal{U}\) is one polyhedral, as shown in Eq.(\ref{eq:D1})-(\ref{eq:D4}). 

\begin{remark}
Clearly, Eq.(\ref{eq:D1}) is similar to Eq.(\ref{eq:gen_status}), regarding the formulation. The failure $\chi_{ij}^t$ is random and the repair $\kappa_{ij}^t$ is given and fixed in Eq.(\ref{eq:D1}), while the start-up $\alpha_{g}^{t}$ and shut down $\beta_{g}^t$ of generators are both schedulable. 
\end{remark}
\begin{remark}
Eq.(\ref{eq:D3}) is equivalent to the following chance constraint 
\begin{equation}\label{eq:D5}
\begin{aligned}
\text{Prob}(\chi_{ij}^{t} \geq 1) \geq \Pi, \forall t,ij
\end{aligned}
\end{equation}

When $\Pi=0$, the uncertainty set \(\mathcal{U}\) reduces to a robust uncertainty set.
\end{remark}
\begin{remark}
For the proposed uncertainty set in Eq.(\ref{eq:D1})-(\ref{eq:D4}), it can be extended to consider the repair operation. When the repair is considered, its repair duration can be formulated as follows
\begin{equation}\label{eq:D6}
\kappa_{ij}^{t+\text{RT}_{ij}} = \chi_{ij}^{t}, \forall t\in\{1,2,...,T-RT_{ij}\},ij
\end{equation}

The uncertainty set with repair is denoted by $\mathcal{U}^{'}$, i.e., Eq.(\ref{eq:D1}), Eq.(\ref{eq:D3}), Eq.(\ref{eq:D4}) and Eq.(\ref{eq:D6}).
\end{remark}

\begin{remark}
To depict the $N-1-1$ contingency reliability criterion, the following constraints can be introduced with $K=2$
\begin{equation}\label{eq:D7}
\sum_{ij \in \mathcal{E}} \chi_{ij}^{t} \leq 1, \forall t
\end{equation}
Eq.(\ref{eq:D7}) indicate at most one line can be failed within each time slot and at most 2 lines are off-line simultaneously.  
\end{remark}

\section{Solution Methods}
The formulated problem (\ref{eq:RO}) is a standard two-stage RO problem. In problem (\ref{eq:RO}), \(\max_{\boldsymbol{\xi} \in \mathcal{U}}[\mathcal{Q}(\textbf{x})]\) is a max-min problem, and the inner problem is a linear programming problem. With Lagrange duality, this max-min problem is reformulated to the following maximal optimization problem 
\begin{equation}\label{eq:dual}
\begin{aligned}
&\mathcal{Q}^{'}(\textbf{x}) = \max_{\boldsymbol{\xi}\in \mathcal{U}, \boldsymbol{\nu}} \boldsymbol{\nu}^{\text{T}}(\textbf{h}-\textbf{G}\textbf{x}-\textbf{M}\boldsymbol{\xi})  \\
& \text{s.t.} \qquad \textbf{E}^{\text{T}}\boldsymbol{\nu} = \textbf{d}, \boldsymbol{\nu} \geq 0
\end{aligned}
\end{equation}

In Eq.(\ref{eq:dual}), for given $\textbf{x}$, it is a non-convex quadratic optimization programming problem, due to the bi-linearity of $\boldsymbol{\nu}^{\text{T}}\textbf{M}\boldsymbol{\xi}$. Considering $\boldsymbol{\xi}$ is a binary vector, $\boldsymbol{\nu}^{\text{T}}\textbf{M}\boldsymbol{\xi}$ can be exactly reformulated using its McCormick envelop\cite{mitsos2009mccormick,zhao2017distributionally}.

The C\&CG scheme \cite{zeng2013solving} is adopted to solve the reformulated problem $\mathcal{Q}^{'}$ as shown in algorithm 1. As $\boldsymbol{\xi}$ is a binary vector, the stopping criterion is $|\boldsymbol{\xi}^{k+1}-\boldsymbol{\xi}^{k}|\leq 1$. 

As the generation can be curtailed and loads can be shedded in Eq.(\ref{eq:generation_curtailment}) and Eq.(\ref{eq:load_shedding}), respectively, $\mathcal{Q}(\textbf{x})$ is a fully recourse problem when $\textbf{X} \neq \emptyset$. Algorithm 1 can converge within finite iterations\cite{zeng2013solving}.
\begin{algorithm}
\SetAlgoLined
\KwData{$\textbf{X}$, $\textbf{Y}$, $\mathcal{U}$, $\textbf{c}$, $\textbf{d}$}
\KwResult{$\textbf{x}$, $\boldsymbol{\xi}$}
 Set $LB$ = $-\inf$, $UB$ = $\inf$, $k$=0, and $\mathcal{O}=\emptyset$
 \While{Stopping criteria not meet}{
  Solve the following master problem
  \begin{equation*}
    \begin{aligned}
    \text{MP}:& \min_{\textbf{x},\eta,, \textbf{y}^{l}, \forall l \leq k} \textbf{c}^{\text{T}}\textbf{x}+\eta \\
    & \text{s.t.} \quad \textbf{G}\textbf{x} + \textbf{E}\textbf{y}^{l}\geq \textbf{h} - \textbf{M}\boldsymbol{\xi}_{l}^{*}, \forall l \leq k\\
    & \qquad \eta \geq \textbf{d}^{\text{T}}\textbf{y}^{l}, \forall l \in \mathcal{O}\\
    & \qquad \textbf{x} \in \textbf{X}, \eta \in \mathbb{R}, \textbf{y}^{l}\in \textbf{Y}, \forall l \leq k
    \end{aligned}
    \end{equation*}
    
    \eIf{MP is infeasible}{Terminate}{
    Derive an optimal solution \(\textbf{x}^{k+1},\eta^{k+1}, \textbf{y}^{1},... \textbf{y}^{k}\).\\
  Update $LB$ = $\textbf{c}^{\text{T}}\textbf{x}^{k+1}+\eta^{k+1}$.\\
  Solve problem (\ref{eq:dual}) and derive $\boldsymbol{\xi}_{k+1}^{*}$.\\
  Update $UB$ = $\textbf{c}^{\text{T}}\textbf{x}^{k+1}+\mathcal{Q}^{'}(\textbf{x}^{k+1})$.
  
  \eIf{$UB-LB \geq \epsilon$}{
   Return $\textbf{x}^{k+1},\boldsymbol{\xi}_{k+1}^{*}$ and terminate
   }{
   Create variables $\textbf{y}^{k+1}$, and add the following constraints to MP
   \begin{equation*}
    \begin{aligned}
    & \eta \geq \textbf{d}^{\text{T}}\textbf{y}^{k+1} \\
    & \textbf{G}\textbf{x} + \textbf{E}\textbf{y}^{k+1}\geq \textbf{h} - \textbf{M}\boldsymbol{\xi}_{k+1}^{*}\\
    \end{aligned}
    \end{equation*}
    Update $k=k+1$, $\mathcal{O}=\mathcal{O}\cup{k+1}$
  }}
 }
 \caption{Column-and-constraint generation algorithm}
\end{algorithm}

\section{Case Studies}
\subsection{Case Description}
Two cases are analyzed to show the effectiveness of the proposed resilient UC in this section. 
\subsubsection{Case I} 
\begin{figure}[!h]
\includegraphics[width=1\linewidth]{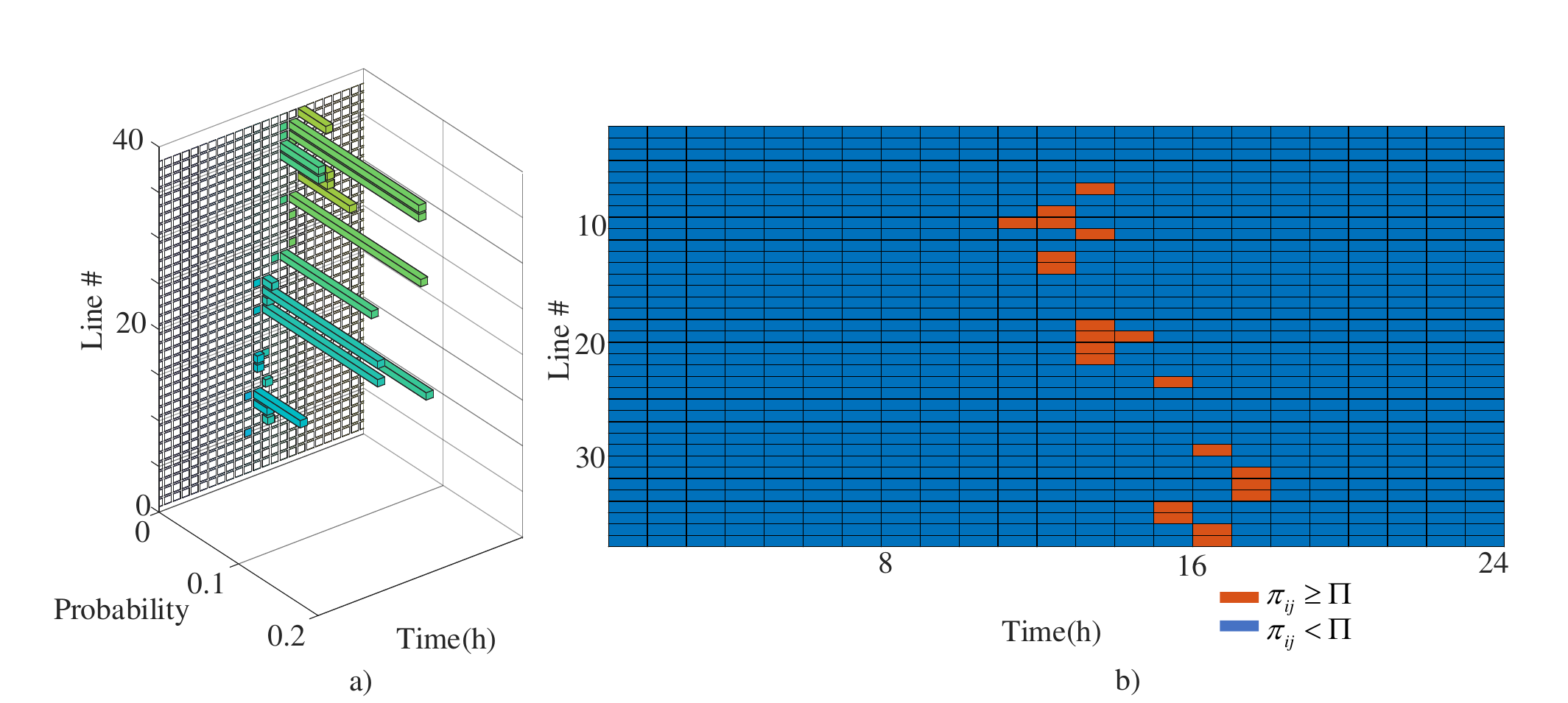}
\caption{Line failure probabilistic of lines during the operating day. a) Line failure probability $\pi_{ij}$. b)$\pi_{ij}>\Pi$.}
\label{fig:linefaileprob}
\end{figure}

A modified IEEE-24 test system under one hurricane is analyzed in this case. The test system is projected to a 150*200 km area
approximately located within (30.52$^{\circ}$N-32.32$^{\circ}$N, 87.68$^{\circ}$W-
89.25$^{\circ}$W). The hurricane is forecast to land within the area (29.31$^{\circ}$N-30.21$^{\circ}$N, 86.64$^{\circ}$W-
90.29$^{\circ}$W) at 8:00 am in the operating day. The hurricane is assumed to have no impact on the power system when it enters the north area of 33.22$^{\circ}$N at 21:00 of the same day. More detailed information on the location and design parameters of towers and conductors are referred to \cite{zhang2019spatial}. The parameters for the generators, loads and transmission lines are referred to \cite{subcommittee1979ieee}. For this hurricane, the failure probability of each line across the hurricane periods is shown in Fig.\ref{fig:linefaileprob}.

\subsubsection{Case II} A two-area IEEE reliability test system-1996\cite{grigg1999ieee} under hurricane is studied to show the effectiveness of proposed UC on interconnected systems. One area has affected by the hurricane, while the other is not affected. The details on other parameters are referred to \cite{zhang2019spatial}. 

Under both cases, $K$ is set to 2, $VOLL$ is set to 4000 \$/MWh, $VOGC$ is set to 1000 \$/MWh, $\Delta t$ is set to 1 hour and $T=24$. Numerical tests were carried out on a laptop with an Intel i7- 4770 CPU and 16 GB of RAM. The optimization problems in algorithm 1 are solved by MILP solver CPLEX\cite{bliek1u2014solving}. To show the effectiveness of the proposed uncertainty set, three scenarios are compared within both cases, as follows
\begin{itemize}
  \item Scenario I, the robust uncertainty set, i.e.,$\Pi=0$.
  \item Scenario II, the proposed uncertainty set, i.e., $\Pi=0.01$.
  \item Scenario III, the repair is considered, i.e., $\mathcal{U}$ is replaced by $\mathcal{U}^{'}$ and $\text{RT}_{ij}=10,\forall ij$.
\end{itemize}

\subsection{Simulation Results of Case I}

The first stage operational cost, load shedding and generation curtailment under the worst line failure scenario are given in Tab.I. As shown in the Tab.I, when the proposed uncertainty set is adopted in scenario II, the first stage operational cost can be reduced by 11.97\% and worst case load shedding can be reduced by 53.87\%, respectively. The operational cost and load shedding can be further reduced when the repair is considered in Scenario III. Under all three cases, there is no generation curtailment when the maximal line failure number is set to 2, i.e., $K=2$. 
\begin{table}[h]
\caption{Simulation results for case I}
\centering
\begin{tabular}{lccc}
\cline{1-3}
\hline
                     & Scenario I & Scenario II & Scenario III \\
First-stage operational cost $\textbf{c}^{\text{T}}\textbf{x}$(\$)         & 81202.67   & 71480.30    & 713515.10    \\
Load shedding(MWh)           & 38.61      & 17.81       & 16.41        \\
Generation curtailment (MWh) & 0          & 0           & 0      \\        \hline
\end{tabular}
\end{table}

\begin{figure}[!h]
\centering
\includegraphics[width=0.7\linewidth]{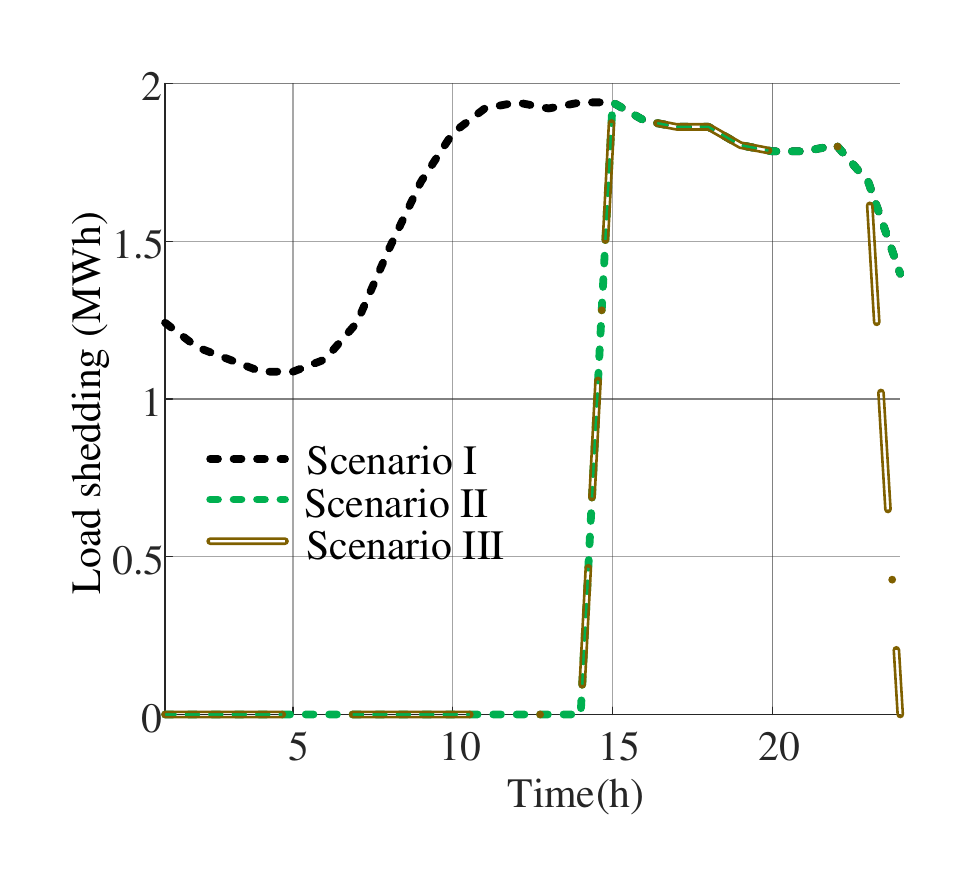}
\caption{Load shedding at bus 14 under worst line failures.}
\label{fig:load_shed}
\end{figure}

\begin{figure}[!h]
\centering
\includegraphics[width=1.05\linewidth]{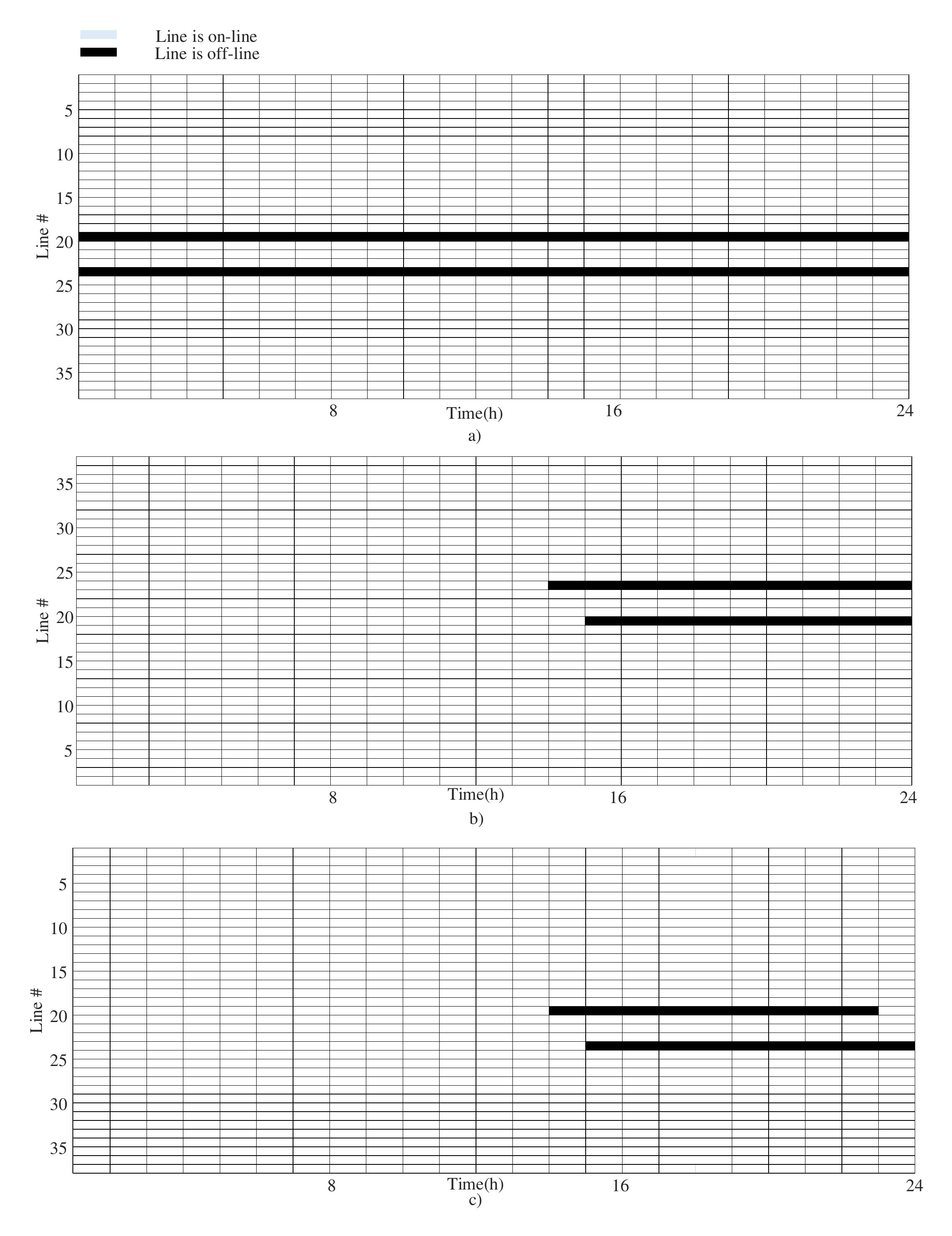}
\caption{Worst line operation scenario for day-ahead market under different scenarios. a) Scenario I, b) Section II, and c) Section III.}
\label{fig:Failure_secnario}
\end{figure}

The load shedding under each scenario is given in Fig.\ref{fig:load_shed} and the worst line failure scenario, i.e., $\textbf{u}^{*}$ in algorithm 1 are depicted in  Fig.\ref{fig:load_shed} and Fig.\ref{fig:Failure_secnario}, respectively. As shown in Fig.\ref{fig:load_shed}, the load at bus 14 is always shedded under all scenario, due to the failure of line 19 (connecting bus 11 and bus 14) and line 23 (connecting bus 14 and bus 16) results in the isolation of bus 14. Considering $\Pi=0.01$, as the increase of failure probability shown in Fig.\ref{fig:linefaileprob} b), line 19 is affected by the hurricane at 14:00 and line 23 is affected by the 15:00, respectively. Line 19 is failed firstly at 14:00, where there is no load shedding at bus 14 at 14:00 as shown in Fig.\ref{fig:load_shed}. Followed by a consecutive failure of line 19 at 15:00, bus 14 is isolated from the system, and the load is shedded to the end the scheduling period under scenario I and scenario II\footnote{This phenomena can be treated as a $N-1-1$ failure event.}. If the lines can be repaired at predefined time, i.e., $\text{RT}_{ij}$, line 19 is on-line at 23:00 as shown in Fig.\ref{fig:Failure_secnario} c); bus 14 is reconnected to the system with sufficient transmission capacity\footnote{$P_{11,14}^{\max} = 500\text{MW}$ and $P_{11,14}^{\min} = -500\text{MW}$} and there is no load shedding at 23:00 under scenario III as shown in Fig.\ref{fig:load_shed}. 

As shown in Eq.(\ref{eq:D4}), parameter $K$ plays an important role regarding the preparedness of the scheduling. The sensitive analysis results with $\Pi=0.01$ are shown in Tab.II. The load shedding amount is not increased linearly with $K$, e.g., $K=3$ and $K=4$ in Tab.II. The first-stage operational cost increases with $K$ before the generation curtailment at $K=7$, where multiple islands are formulated, and the generation should be curtailed to meet the energy balance within each sub-system.

\begin{table}[h]
\caption{Results under different $K$}
\centering
\begin{tabular}{lccc}
\hline
                                & $K$=2      & $K$=3      & $K$=4      \\
\hline
First-stage operational cost (\$) & 71480.30   & 814538.20 & 814538.20 \\
Load shedding (MWh)             & 17.81  & 25.28  & 25.28  \\
Generation curtailment (MWh)    & 0        & 0        & 0        \\
\hline
                                & $K$=5      & $K$=6      & $K$=7      \\
                                & 826920.20 & 826920.20 & 762577.00   \\
                                & 43.08  & 43.08  & 41.59  \\
                                & 0        & 0        & 17.41 \\
\hline
\end{tabular}
\end{table}

\begin{figure}[!h]
\centering
\includegraphics[width=0.8\linewidth]{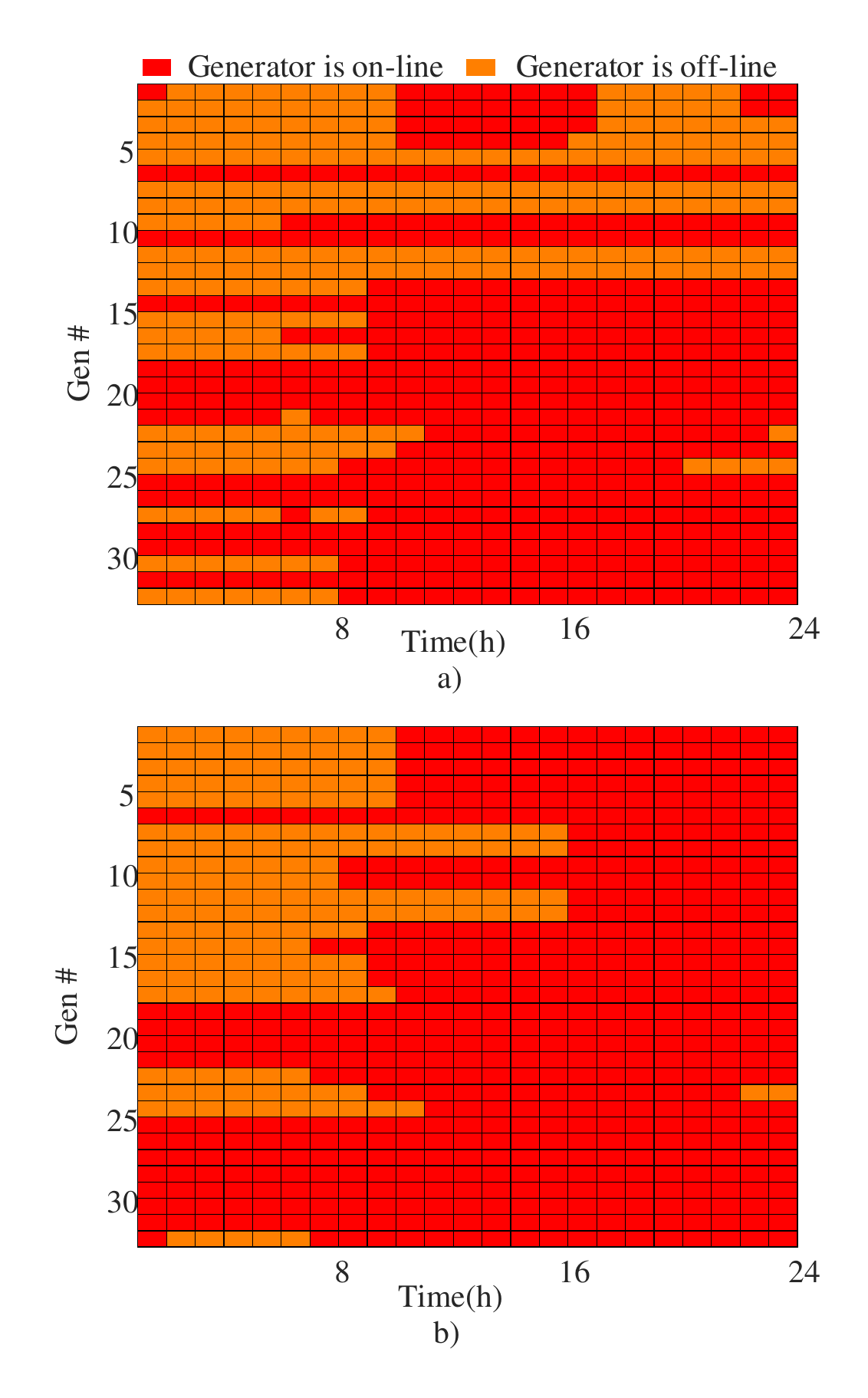}
\caption{Unit operation status under different $K$. a) $K$=2, b) $K$=7.}
\label{fig:gen_status}
\end{figure}

After the clear of the day-ahead market, the unit operation status is shown in Fig.\ref{fig:gen_status}. As shown in Fig.\ref{fig:gen_status}, more generators are on-line continuously with less start-up and shut-down operations, reducing the operational cost in the day-ahead market, as shown in Tab.II. The on-line generators might face curtailment when 7 transmission lines fail in the operating day, as shown in Tab.II when $K=7$. 

The results indicate the proposed uncertainty set can: 1) reduce the conservation of robust uncertainty set, as shown in Tab.I; 2) quantify the impacts of repair on the operational cost and load shedding, as shown in Tab.I; 3) detect worst failure event under given criterion, as shown in Fig.\ref{fig:load_shed}; and  4) provide guidance to generators in face of severe failures, as shown in Tab.II.

\subsection{Simulation Results of Case II}
The simulation results for case II are given in Tab.III. As shown in Tab.III, the interconnected two area system is more resilient regarding the hurricane, as the operational cost fluctuates within a smaller range, in comparison with the result in Tab.I. What is more, the same worst line failures have been found in the area affected by the hurricane, indicated by the same load shedding amount in both case I and case II.

\begin{table}[h]
\caption{Results of case II}
\centering
\begin{tabular}{lccc}
\cline{1-3}
\hline
                             & Scenario I & Scenario II & Scenario III \\
First-stage operational cost $\textbf{c}^{\text{T}}\textbf{x}$(\$)         & 1402146.62   & 1402670.31    & 1399924.26    \\
Load shedding(MWh)           & 38.61      & 17.81       & 16.41        \\
Generation curtailment (MWh) & 0          & 0           & 0      \\        \hline
\end{tabular}
\end{table}

\begin{figure}[!h]
\centering
\includegraphics[width=0.8\linewidth]{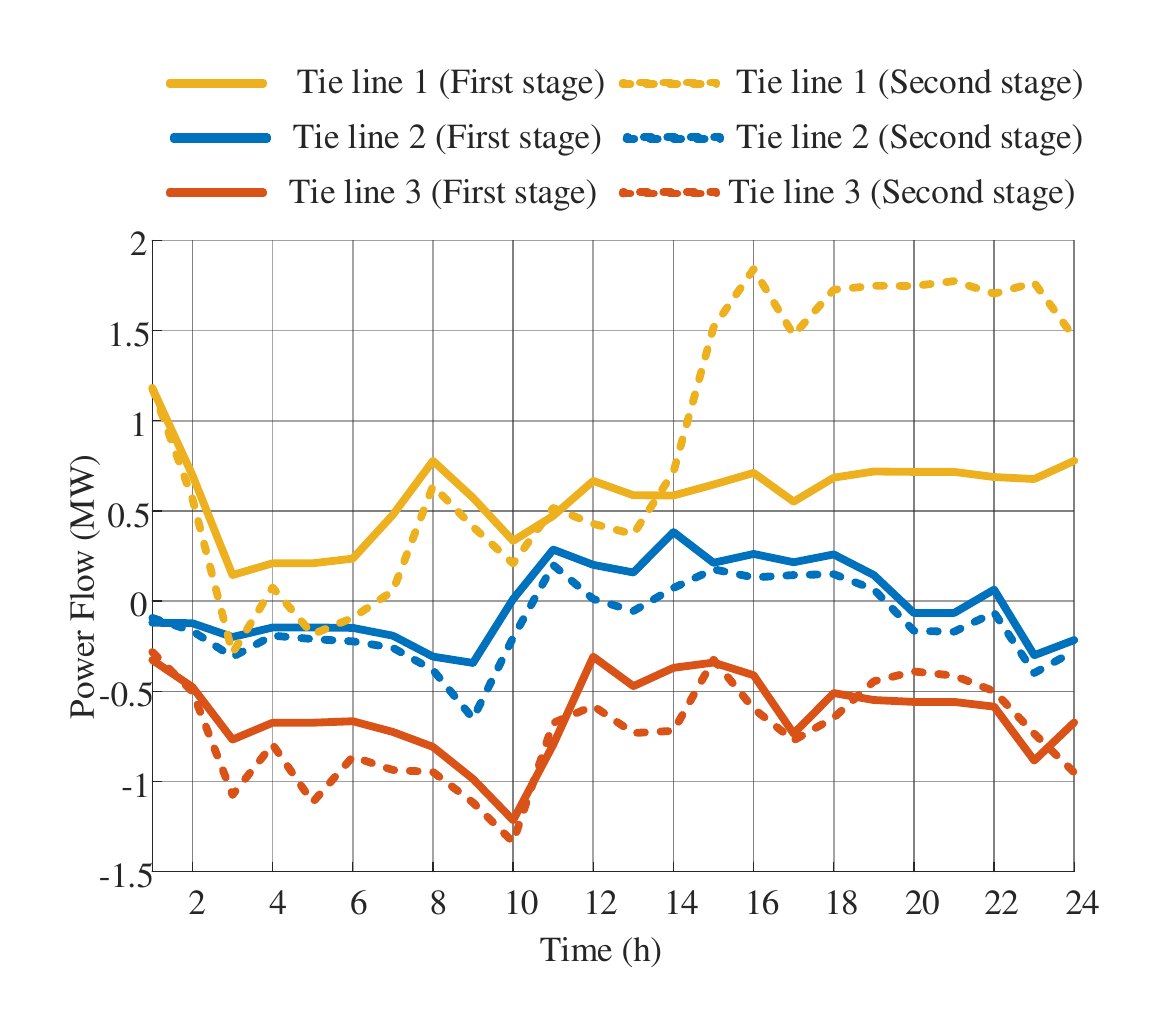}
\caption{Power flow on tie-lines in the first stage and second stage.}
\label{fig:tie_lins}
\end{figure}

The tie-line power flows in the first stage and worst second-stage line failures are shown in Fig.\ref{fig:tie_lins}. The power flow of tie-line 1 will be increased from the affected area to the unaffected area, to absorb the excess generation due to load shedding. It indicates the impact propagation by hurricanes within the interconnected power systems.
 
\section{Conclusion}
In this paper, a resilient UC is proposed for the day-ahead market towards hurricanes. The impacts of hurricanes on transmission lines are formulated as Bernoulli distribution with time-varying parameters. Integrating the line failure probability and introducing the repair and failure operation into a robust uncertainty set with chance constraints, a two-stage robust unit commitment problem is formulated to enhance the resilience of power systems. The problem is solved using the column-and-constraint generation scheme. Simulation is performed on the modified IEEE-24 system and two-area IEEE reliability test system-1996. Results indicate the proposed uncertainty set can reduce the conservation of robust line failure uncertainty set via considering the time and space varying line failures probability, detect the worst line failure event, assess the effectiveness of the repair, and provide guidance for generators and tie-lines among interconnected power systems.

\bibliographystyle{IEEEtran}
\bibliography{references}

\begin{thebibliography}{10}
\providecommand{\url}[1]{#1}
\csname url@samestyle\endcsname
\providecommand{\newblock}{\relax}
\providecommand{\bibinfo}[2]{#2}
\providecommand{\BIBentrySTDinterwordspacing}{\spaceskip=0pt\relax}
\providecommand{\BIBentryALTinterwordstretchfactor}{4}
\providecommand{\BIBentryALTinterwordspacing}{\spaceskip=\fontdimen2\font plus
\BIBentryALTinterwordstretchfactor\fontdimen3\font minus
  \fontdimen4\font\relax}
\providecommand{\BIBforeignlanguage}[2]{{%
\expandafter\ifx\csname l@#1\endcsname\relax
\typeout{** WARNING: IEEEtran.bst: No hyphenation pattern has been}%
\typeout{** loaded for the language `#1'. Using the pattern for}%
\typeout{** the default language instead.}%
\else
\language=\csname l@#1\endcsname
\fi
#2}}
\providecommand{\BIBdecl}{\relax}
\BIBdecl

\bibitem{zhang2019spatial}
H.~Zhang, L.~Cheng, S.~Yao, T.~Zhao, and P.~Wang, ``Spatial-temporal
  reliability and damage assessment of transmission networks under
  hurricanes,'' \emph{IEEE Transactions on Smart Grid}, 2019.

\bibitem{wang2015research}
Y.~Wang, C.~Chen, J.~Wang, and R.~Baldick, ``Research on resilience of power
  systems under natural disasters—a review,'' \emph{IEEE Transactions on
  Power Systems}, vol.~31, no.~2, pp. 1604--1613, 2015.

\bibitem{zhao2017distributionally}
C.~Zhao and R.~Jiang, ``Distributionally robust contingency-constrained unit
  commitment,'' \emph{IEEE Transactions on Power Systems}, vol.~33, no.~1, pp.
  94--102, 2017.

\bibitem{ela2014evolution}
E.~Ela, M.~Milligan, A.~Bloom, A.~Botterud, A.~Townsend, and T.~Levin,
  ``Evolution of wholesale electricity market design with increasing levels of
  renewable generation,'' National Renewable Energy Lab.(NREL), Golden, CO
  (United States), Tech. Rep., 2014.

\bibitem{fernandez2016probabilistic}
R.~Fern{\'a}ndez-Blanco, Y.~Dvorkin, and M.~A. Ortega-Vazquez, ``Probabilistic
  security-constrained unit commitment with generation and transmission
  contingencies,'' \emph{IEEE Transactions on Power Systems}, vol.~32, no.~1,
  pp. 228--239, 2016.

\bibitem{zheng2014stochastic}
Q.~P. Zheng, J.~Wang, and A.~L. Liu, ``Stochastic optimization for unit
  commitment—a review,'' \emph{IEEE Transactions on Power Systems}, vol.~30,
  no.~4, pp. 1913--1924, 2014.

\bibitem{Dimitris2019Resilient}
D.~N.~Trakas and N.~D.~Hatziargyriou, ``Resilience constrained day-ahead unit
  commitment under extreme weather events,'' \emph{IEEE Transactions on Power
  Systems}, 2019.

\bibitem{wang2016resilience}
C.~Wang, Y.~Hou, F.~Qiu, S.~Lei, and K.~Liu, ``Resilience enhancement with
  sequentially proactive operation strategies,'' \emph{IEEE Transactions on
  Power Systems}, vol.~32, no.~4, pp. 2847--2857, 2016.

\bibitem{wang2018resilience}
Y.~Wang, L.~Huang, M.~Shahidehpour, L.~L. Lai, H.~Yuan, and F.~Y. Xu,
  ``Resilience-constrained hourly unit commitment in electricity grids,''
  \emph{IEEE Transactions on Power Systems}, vol.~33, no.~5, pp. 5604--5614,
  2018.

\bibitem{wang2013two}
Q.~Wang, J.-P. Watson, and Y.~Guan, ``Two-stage robust optimization for nk
  contingency-constrained unit commitment,'' \emph{IEEE Transactions on Power
  Systems}, vol.~28, no.~3, pp. 2366--2375, 2013.

\bibitem{guo2016contingency}
Z.~Guo, R.~L.-Y. Chen, N.~Fan, and J.-P. Watson, ``Contingency-constrained unit
  commitment with intervening time for system adjustments,'' \emph{IEEE
  Transactions on Power Systems}, vol.~32, no.~4, pp. 3049--3059, 2016.

\bibitem{hu2015robust}
B.~Hu and L.~Wu, ``Robust scuc considering continuous/discrete uncertainties
  and quick-start units: A two-stage robust optimization with mixed-integer
  recourse,'' \emph{IEEE Transactions on Power Systems}, vol.~31, no.~2, pp.
  1407--1419, 2015.

\bibitem{wang2016robust}
C.~Wang, W.~Wei, J.~Wang, F.~Liu, F.~Qiu, C.~M. Correa-Posada, and S.~Mei,
  ``Robust defense strategy for gas--electric systems against malicious
  attacks,'' \emph{IEEE Transactions on Power Systems}, vol.~32, no.~4, pp.
  2953--2965, 2016.

\bibitem{babaei2018distributionally}
S.~Babaei, R.~Jiang, and C.~Zhao, ``Distributionally robust distribution
  network configuration under random contingency,'' \emph{arXiv preprint
  arXiv:1808.09038}, 2018.

\bibitem{hurricane_report}
``Hurricane harvey event analysis report,'' North American Reliability
  Corporation.(NERC), Washington, DC (United States), Tech. Rep., 2018.

\bibitem{subcommittee1979ieee}
P.~M. Subcommittee, ``Ieee reliability test system,'' \emph{IEEE Transactions
  on power apparatus and systems}, no.~6, pp. 2047--2054, 1979.

\bibitem{van2015transmission}
P.~Van~Hentenryck and C.~Coffrin, ``Transmission system repair and
  restoration,'' \emph{Mathematical Programming}, vol. 151, no.~1, pp.
  347--373, 2015.

\bibitem{eu_market_report}
``Overview of european electricity markets,'' European Union, Tech. Rep., 2016.

\bibitem{noauthor_despite_nodate}
``Despite customer outages, wholesale electric markets operated during
  {Hurricane} {Sandy} - {Today} in {Energy} - {U}.{S}. {Energy} {Information}
  {Administration} ({EIA}).''

\bibitem{sandi_report}
``Hurricane sandy event analysis report,'' North American Reliability
  Corporation.(NERC), Washington, DC (United States), Tech. Rep., 2014.

\bibitem{yao2018transportable}
S.~Yao, P.~Wang, and T.~Zhao, ``Transportable energy storage for more resilient
  distribution systems with multiple microgrids,'' \emph{IEEE Transactions on
  Smart Grid}, vol.~10, no.~3, pp. 3331--3341, 2018.

\bibitem{singhal_hedman_2019}
N.~Singhal, J.~Kwon, and K.~Hedman, ``Generator contingency modeling in
  electric energy markets: Derivation of prices via duality theory,''
  \emph{arXiv preprint arXiv:1910.02323}, 2019.

\bibitem{Xiao2006c}
F.~Xiao, J.~D. McCalley, Y.~Ou, J.~Adams, and S.~Myers, ``{Contingency
  Probability Estimation Using Weather and Geographical Data for On-line
  Security Assessment},'' in \emph{2006 Int. Conf. Probabilistic Methods Appl.
  to Power Syst.}\hskip 1em plus 0.5em minus 0.4em\relax IEEE, 2006, pp. 1--7.

\bibitem{Liu2011f}
Y.~Liu and C.~Singh, ``{A Methodology for Evaluation of Hurricane Impact on
  Composite Power System Reliability},'' \emph{IEEE Trans. Power Syst.},
  vol.~26, no.~1, pp. 145--152, 2011.

\bibitem{morales2012tight}
G.~Morales-Espa{\~n}a, J.~M. Latorre, and A.~Ramos, ``Tight and compact milp
  formulation of start-up and shut-down ramping in unit commitment,''
  \emph{IEEE Transactions on Power Systems}, vol.~28, no.~2, pp. 1288--1296,
  2012.

\bibitem{li2015flexible}
N.~Li, C.~Uckun, E.~M. Constantinescu, J.~R. Birge, K.~W. Hedman, and
  A.~Botterud, ``Flexible operation of batteries in power system scheduling
  with renewable energy,'' \emph{IEEE Transactions on Sustainable Energy},
  vol.~7, no.~2, pp. 685--696, 2015.

\bibitem{mitsos2009mccormick}
A.~Mitsos, B.~Chachuat, and P.~I. Barton, ``Mccormick-based relaxations of
  algorithms,'' \emph{SIAM Journal on Optimization}, vol.~20, no.~2, pp.
  573--601, 2009.

\bibitem{zeng2013solving}
B.~Zeng and L.~Zhao, ``Solving two-stage robust optimization problems using a
  column-and-constraint generation method,'' \emph{Operations Research
  Letters}, vol.~41, no.~5, pp. 457--461, 2013.

\bibitem{grigg1999ieee}
C.~Grigg, P.~Wong, P.~Albrecht, R.~Allan, M.~Bhavaraju, R.~Billinton, Q.~Chen,
  C.~Fong, S.~Haddad, S.~Kuruganty \emph{et~al.}, ``The ieee reliability test
  system-1996. a report prepared by the reliability test system task force of
  the application of probability methods subcommittee,'' \emph{IEEE
  Transactions on power systems}, vol.~14, no.~3, pp. 1010--1020, 1999.

\bibitem{bliek1u2014solving}
C.~Bliek1{\'u}, P.~Bonami, and A.~Lodi, ``Solving mixed-integer quadratic
  programming problems with ibm-cplex: a progress report,'' in
  \emph{Proceedings of the twenty-sixth RAMP symposium}, 2014, pp. 16--17.

\end{thebibliography}

\end{document}